# STRANGER THINGS ABOUT FORCING WITHOUT AC

MARTIN GOLDSTERN AND LUKAS DANIEL KLAUSNER

ABSTRACT. Typically, set theorists reason about forcing constructions in the context of ZFC. We show that without AC, several simple properties of forcing posets fail to hold, one of which answers Miller's question from [Mil08].

## 1. Fin$(X, 2)$ AND CARDINAL COLLAPSE

Miller [Mil08, p. 314] posed the question whether forcing with Fin$(X, 2)$ in ZF can make two sets $A$ and $B$ of different cardinality in the ground model have the same cardinality in the generic extension. We find that "collapses" are possible in the sense that non-equipotent sets may become equipotent after forcing with Fin$(X, 2)$, but not if both sets are well-ordered.

**Definition 1.1.** We write $A \approx B$ to abbreviate "there is a bijection from $A$ onto $B$".

**Example 1.2.** *Let $(A_n)_{n<\omega}$ be a countable sequence of pairs of socks, i. e. a sequence of pairwise disjoint two-element sets which does not have a choice function. Let $A := \bigcup_{n<\omega} A_n$ and $\mathbb{P} := \mathrm{Fin}(A \times \omega, 2)$. Then we have that $V \vDash A \not\approx \omega$ (since there is no choice function), but $V^{\mathbb{P}} \vDash A \approx \omega$.*

*Similarly, with $X := A \times \omega$ we get that Fin$(X, 2)$ forces $X \approx \omega$, while $X \not\approx \omega$ in the ground model.*

*Proof.* Let $D_n := \{p \in \mathbb{P} \mid \exists\, m < \omega \colon p[A_n \times \{m\}] = \{0, 1\}\}$. For any $n$, $D_n$ is dense in $\mathbb{P}$: Given $p \in \mathbb{P}$, let $m$ be minimal with $p\restriction_{A_n \times \{m\}} = \varnothing$. Letting $A_n = \{a, b\}$ (choosing once), define $q := p \cup \{((a, m), 0), ((b, m), 1)\}$. Then $p \geq q \in D_n$.

The generic $g$ then satisfies that for any $n < \omega$, there is an $m < \omega$ such that $g[A_n \times \{m\}] = 2$. We now choose for each $n < \omega$ the minimal $m_n < \omega$ with that property, and thus have a well-order on $A$ in the generic extension (by letting $A_n = \{a, b\}$ be ordered by $a < b \Leftrightarrow g(a, m_n) < g(b, m_n)$). □

**Lemma 1.3.** *Let $A$ be an antichain in* Fin$(X, 2)$ *consisting only of conditions with domains of size $k$. Then $A$ is finite, and $|A|$ is even bounded by $2^k$.*

*Proof.* It suffices to show that every finite subset $A' \subseteq A$ has at most $2^k$ elements.

So let $A' \subseteq A$ be finite with $n$ elements. Let $D := \bigcup_{a \in A'} \mathrm{dom}(a)$ and $d := |D|$. Consider the set $2^D$; for each $a \in A'$, there are exactly $2^{d-k}$ many possible extensions of $a$ in $2^D$, and since $A'$ is an antichain, they must all be pairwise different.

2010 *Mathematics Subject Classification.* Primary 03E25; Secondary 03E40.
Both authors were supported by the Austrian Science Fund (FWF) project P29575 "Forcing Methods: Creatures, Products and Iterations". We are grateful to the anonymous referee for suggesting improvements, and in particular for asking a question that led to Remark 1.5.

Hence there are $n \cdot 2^{d-k}$ many different such elements of $2^D$, thus $n \cdot 2^{d-k} \leq 2^d$, which concludes the proof.[1] $\square$

**Theorem 1.4.** *Let $\kappa < \lambda$ be well-ordered cardinals. Then there is no $X$ such that* $\mathrm{Fin}(X, 2)$ *forces $\kappa \approx \lambda$.*

*Proof.* Let $\dot{f}$ be a $\mathrm{Fin}(X, 2)$-name and let $p \in \mathrm{Fin}(X, 2)$ be such that $p \Vdash$ "$\dot{f} \colon \kappa \to \lambda, \dot{f}$ is onto". For each $\alpha \in \kappa$ and $k < \omega$, define

$$A_{\alpha,k} := \{\beta \in \lambda \mid \exists\, q \leq p \colon |\mathrm{dom}(q)| = k, q \Vdash \dot{f}(\alpha) = \beta\}.$$

We claim that all such $A_{\alpha,k}$ are finite, even bounded by $2^k$: Assume otherwise, i. e. let $\beta_1, \ldots, \beta_{2^k+1}$ be different elements of $A_{\alpha,k}$. Then there are witnesses $q_1, \ldots, q_{2^k+1}$ to that, which necessarily are incompatible and thus form a $(2^k+1)$-sized antichain of conditions with domains of size $k$, giving an obvious contradiction to Lemma 1.3.

From here on, we work in the ground model. By the above, we have that

$$\lambda = \bigcup_{\alpha \in \kappa} \bigcup_{k < \omega} A_{\alpha,k}.$$

Using the bijection $\varphi \colon \kappa \times \omega \to \kappa \colon (\alpha, k) \mapsto \omega \cdot \alpha + k$ to arrange the sets linearly, we have that

$$\lambda = \bigcup_{\gamma \in \kappa} B_\gamma$$

with finite $B_\gamma$. As the $B_\gamma$ are subsets of $\lambda$ (thus well-ordered) as well as finite, we can embed $\bigcup_{\gamma \in \kappa} B_\gamma$ into $\kappa \times \omega$ (by embedding $B_\gamma$ into $\{\gamma\} \times \omega$) and hence into $\kappa$ (again using $\varphi$), arriving at a contradiction. $\square$

We remark that a similar argument shows that the statement "$\mathrm{cof}(\alpha) = \beta$" is absolute for any pair of ordinals $(\alpha, \beta)$.

**Remark 1.5.** Can the forcing $\mathrm{Fin}(X)$ add a bijection between two sets $A$, $B$ for which there is no such bijection in the ground model? In the following (rough) analysis we will write $\mathrm{WO}(X)$ to abbreviate "$X$ can be well-ordered".

- Assume $\neg\,\mathrm{WO}(X)$.
  - Example 1.2 shows that the answer can be "yes" if one of $\neg\,\mathrm{WO}(A)$, $\neg\,\mathrm{WO}(B)$ holds. A minor variation shows that this may also be possible if both $\neg\,\mathrm{WO}(A)$ and $\neg\,\mathrm{WO}(B)$ hold.
  - If both $\mathrm{WO}(A)$ and $\mathrm{WO}(B)$ hold, then the answer is "no", as shown in Theorem 1.4.
- Now assume $\mathrm{WO}(X)$.
  - If both $\mathrm{WO}(A)$ and $\mathrm{WO}(B)$ hold, then the answer is "no" once again.
  - Asaf Karagila [Kar18] pointed out that the answer is also "no" if $\mathrm{WO}(A)$ and $\neg\,\mathrm{WO}(B)$ hold.
  - While we suspect that the answer is also "no" if both $\neg\,\mathrm{WO}(A)$ and $\neg\,\mathrm{WO}(B)$ hold, we cannot prove this yet.

[1] We wish to thank Martin Ziegler for pointing out this simple proof for this (optimal) bound for the size of $A$.

## 2. $\sigma$-Closedness

**Definition 2.1.** A forcing poset $\mathbb{P}$ is $\sigma$-closed if for all descending sequences of conditions $\langle p_n \mid n < \omega \rangle$ in $\mathbb{P}$, there is a $q \in \mathbb{P}$ with $q \leq p_i$ for all $i < \omega$.

It is well-known that in ZFC, $\sigma$-closed forcing posets have a number of "nice" properties; the following examples (which are probably part of set-theoretic folklore) show that this is not true in ZF alone.

**Proposition 2.2.** *Let $X$ be a Dedekind-finite infinite set (i. e. there is no $\omega$-sequence within $X$, or equivalently, there is no injective function from $X$ onto a proper subset of $X$). Then the forcing poset $\mathbb{P} := \mathrm{Fin}_{\mathrm{inj}}(\omega, X)$ of partial finite injective functions from $\omega$ to $X$ is $\sigma$-closed.*[2]

*Proof.* It is well-known that $P$ is Dedekind-finite, hence trivially $\sigma$-closed. For completeness' sake, we give an explicit proof:

Let $\langle p_n \mid n < \omega \rangle$ be a descending sequence of conditions in $\mathbb{P}$. We claim that the sequence must be eventually constant; if that is the case and $\langle p_n \mid n < \omega \rangle$ is eventually constant beginning with $p_k = p_{k+1} = p_{k+2} = \ldots$, let $q := p_k$.

To prove the claim, assume that the sequence is not eventually constant, i. e. there is a strictly monotone sequence of integers $\langle \ell(n) \mid n < \omega \rangle$ such that $p_{\ell(0)} > p_{\ell(1)} > p_{\ell(2)} > \ldots$. For all $n < \omega$, let $s(n)$ be defined as the first new element in the domain of $p_{\ell(n)}$, i. e.

$$s(n) := \min \left( \mathrm{dom}(p_{\ell(n)}) \setminus \bigcup_{k<n} \mathrm{dom}(p_{\ell(k)}) \right).$$

But then $x_n := p_{\ell(n)}(s(n))$ is an $\omega$-sequence within $X$, contradicting the assumption. $\square$

Using this, we can define forcing posets which are $\sigma$-closed but add new reals and/or collapse cardinals. Let $\mathrm{Fin}_{\pi_1\text{-inj}}(A, B \times C)$ denote the poset of partial finite functions from $A$ to $B \times C$ which are $\pi_1$-injective (where $\pi_1$ is the projection onto the first coordinate), i. e. functions $f$ for which $\pi_1 \circ f$ (from $A$ to $B$) is injective.

**Example 2.3.** *Let $X$ be a Dedekind-finite infinite set (i. e. there is no $\omega$-sequence within $X$, or equivalently, there is no injective function from $X$ onto a true subset of $X$). Then the forcing posets $\mathbb{P}_1 := \mathrm{Fin}_{\pi_1\text{-inj}}(\omega, X \times 2)$ and $\mathbb{P}_2 := \mathrm{Fin}_{\pi_1\text{-inj}}(\omega, X \times \omega_1)$ are $\sigma$-closed, but $\mathbb{P}_1$ adds a new real and $\mathbb{P}_2$ collapses $\omega_1^V$ to $\omega^V$.*[3]

*Proof.* First, for $\mathbb{P}_1$, assume we are given a descending sequence of conditions $\langle (p_n^X, p_n^2) \mid n < \omega \rangle$ (writing $(p^X, p^2)$ for $(\pi_1 \circ p, \pi_2 \circ p)$). By Proposition 2.2, the sequence $\langle p_n^X \mid n < \omega \rangle$ must be eventually constant, hence so must their domains, and hence so must the sequence $\langle p_n^2 \mid n < \omega \rangle$. Therefore, $\mathbb{P}_1$ is $\sigma$-closed. However, due to the second components of the conditions, forcing with $\mathbb{P}_1$ adds a Cohen real.

[2] Alternatively, we could consider the poset $\mathbb{P} := \mathrm{Fin\,Seq}_{\mathrm{inj}}(X)$ of finite sequences of different elements in $X$: If $X$ is $\sigma$-closed, then so is $\mathbb{P}$.

[3] Of course, $\mathbb{P}_2$ then *also* adds a new real (namely the well-order of $\omega_1^V$).

For $\mathbb{P}_2$ the same argument shows $\sigma$-closedness. Here, the second components of the conditions ensure that forcing with $\mathbb{P}_2$ adds a function from $\omega^V$ onto $\omega_1^V$ and hence collapses $\omega_1^V$. $\square$

Sufficiently disenchanted by this result, we would like to introduce a property of forcing posets which

- is equivalent to $\sigma$-closedness under DC and
- implies that the forcing poset adds no new reals even in ZF.

**Definition 2.4.** A family $\langle P_n \mid n < \omega\rangle$ of disjoint subsets of a forcing poset $\mathbb{P}$ is a *pyramid* (in $\mathbb{P}$) if for all $n < \omega$ and all $p \in P_n$, there is a $k > n$ and a $q \in P_k$ such that $q \leq p$.

A forcing poset $\mathbb{P}$ is *capstone-closed* if for each pyramid $\langle P_n \mid n < \omega\rangle$ in $\mathbb{P}$, there is a *capstone* $q \in \mathbb{P}$ such that for all $k < \omega$, there is a $p \in \bigcup_{n\geq k} P_k$ with $q \leq p$.

**Lemma 2.5.** *Assume* DC. *If a forcing poset $\mathbb{P}$ is $\sigma$-closed, then it is also capstone-closed.*

*Proof.* Let $\langle P_n \mid n < \omega\rangle$ be a pyramid in $\mathbb{P}$. Use DC and the defining property of a pyramid to find a strictly increasing sequence $\langle k_n \mid n < \omega\rangle$ and a descending sequence of conditions $\langle p_n \mid n < \omega\rangle$ in $\mathbb{P}$ such that $p_n \in P_{k_n}$ for all $n < \omega$. By the $\sigma$-closedness of $\mathbb{P}$, there must be some $q \in \mathbb{P}$ such that $q \leq p_n$ for all $n < \omega$. This $q$ then is the capstone of the pyramid $\langle P_n \mid n < \omega\rangle$ as witnessed by the $p_n$. $\square$

The converse holds even in ZF:

**Lemma 2.6.** *If a forcing poset $\mathbb{P}$ is capstone-closed, then it is also $\sigma$-closed.*

*Proof.* Let $\langle p_n \mid n < \omega\rangle$ be a descending sequence of conditions in $\mathbb{P}$. Without loss of generality, assume it is strictly descending. The capstone $q$ of the pyramid $\langle \{p_n\} \mid n < \omega\rangle$ fulfils $q \leq p_n$ for all $n < \omega$. $\square$

Finally, we show that capstone-closedness ensures that forcing posets behave nicely even in ZF.

**Theorem 2.7.** *If a forcing poset $\mathbb{P}$ is capstone-closed, it adds neither new reals nor new sequences of ordinals.*

*Proof.* Assume towards a contradiction that there is a condition $p^* \in \mathbb{P}$ and a $\mathbb{P}$-name $\dot{x}$ such that $p^* \Vdash$ "$\dot{x}\colon \omega \to \mathrm{Ord}$ is new".

Let $P_n$ be the set of all conditions below $p^*$ which decide $\dot{x}\restriction_n$, but not $\dot{x}\restriction_{n+1}$. Then $\langle P_n \mid n < \omega\rangle$ is a disjoint covering of the conditions below $p^*$: It is clear that no condition can be in more than one $P_n$, and if $p$ were in none of the $P_n$, it would decide all of $\dot{x}$ – but then it would force that $\dot{x}$ is not a new real, which is a contradiction. Moreover, $\langle P_n \mid n < \omega\rangle$ is a pyramid: Given $n < \omega$ and $p \in P_n$, there is some $q \leq p$ deciding $\dot{x}\restriction_{n+1}$; since $q$ cannot decide all of $\dot{x}$, there must be some $k > n$ such that $q \in P_k$.

Now let $q$ be the capstone of the pyramid $\langle P_n \mid n < \omega\rangle$. For each $n < \omega$, there is exactly one $s_n \in \mathrm{Ord}^n$ and a condition $p$ such that $q \leq p$ and such that $p \Vdash s_n \subseteq$ $\dot{x}$. These $s_n$ are necessarily pairwise compatible, hence $s^* := \bigcup_{n<\omega} s_n$ exists and $q \Vdash \dot{x} = s^*$ – but $s^*$ is in the ground model, which is a contradiction. □

## 3. Finite Chain Condition

**Definition 3.1.** A forcing poset $\mathbb{P}$ has the $\theta$-chain condition (is $\theta$-cc) if every family of incompatible elements of $\mathbb{P}$ has cardinality $< \theta$.

It is well known that (in ZFC) every forcing poset that is $\aleph_0$-cc must also be $n$-cc for some $n < \omega$ (see [Kun83, Exercises for Chapter (VII) (F1)]). We show that this is not true in ZF.

**Example 3.2.** *Let $X$ be infinite without a function from $X$ onto $\omega$ (for example, let $X$ be an amorphous set); this is consistent with* ZF. *Let $\mathbb{P} := \mathrm{Fin}(X, 2)$. For any $n$, there is an $n$-element subset $E$ of $X$, and the set of functions from $E$ into 2 is an antichain of size $2^n$ in $\mathbb{P}$; however, there is no infinite antichain.*

*Proof.* Let us assume towards a contradiction that $C$ is an infinite antichain, and let $C_n$ be the set of elements of $C$ of size at most $n$. Let $X_n := \bigcup_{c \in C_n} \mathrm{dom}(c)$. By Lemma 1.3, each $C_n$ is finite, and so is each $X_n$.

Now define $f : \mathbb{P} \to \omega$ by defining $f{\restriction}_{\bigcup_{n<\omega} X_n}$ as $f(x) := \min\{n \mid x \in X_n\}$ and letting $f{\restriction}_{\mathbb{P} \setminus \bigcup_{n<\omega} X_n} \equiv 0$. Then $f$ is a map from $\mathbb{P}$ onto an infinite subset of $\omega$, so there is some $g : \omega \to \omega$ such that $g \circ f$ is onto, which is a contradiction. □

## References

[Kar18] Asaf Karagila, *Do choice principles in all generic extensions imply AC in V?*, answer to a MathOverflow question, 2018, MathOverflow: 307269.

[Kun83] Kenneth Kunen, *Set Theory: An Introduction to Independence Proofs*, Studies in Logic and the Foundations of Mathematics, vol. 102, North-Holland, Amsterdam, 1983.

[Mil08] Arnold W. Miller, *Long Borel Hierarchies*, MLQ Math. Log. Q. **54** (2008), no. 3, 307–322, DOI: 10.1002/MALQ.200710044, arXiv: 0704.3998 [math.LO].

Institute of Discrete Mathematics and Geometry, TU Wien, Wiedner Hauptstrasse 8–10/104, 1040 Wien, Austria

*Email address*: goldstern@tuwien.ac.at

*URL*: http://www.tuwien.ac.at/goldstern/

Institute of Discrete Mathematics and Geometry, TU Wien, Wiedner Hauptstrasse 8–10/104, 1040 Wien, Austria

*Email address*: mail@l17r.eu

*URL*: https://l17r.eu